\documentclass[12pt,leqno]{article}
\usepackage{amsmath,amssymb}
\begin{document}
\newtheorem{theorem}{Theorem}[section]
%\label %\ref
\newtheorem{prop}{Proposition}[section]
\newtheorem{defin}{Definition}[section]
\newtheorem{rem}{Remark}[section]
\newtheorem{example}{Example}[section]
\newtheorem{corollary}{Corollary}[section]
%\label{} %\ref{}
\title{K\"ahler-Nijenhuis Manifolds}
\author{{\small by}\vspace{2mm}\\Izu Vaisman}
\date{}
\maketitle
{\def\thefootnote{*}\footnotetext[1]%
{{\it 2000 Mathematics Subject Classification: 53C55, 53D17} .
Poisson-Nijenhuis structure. K\"ahler-Nijenhuis manifold.}}
\begin{center} \begin{minipage}{12cm}
A{\footnotesize BSTRACT. A K\"ahler-Nijenhuis manifold is a
K\"ahler manifold $M$, with metric $g$, complex structure $J$ and
K\"ahler form $\Omega$, endowed with a Nijenhuis tensor field $A$
that is compatible with the Poisson structure defined by $\Omega$
in the sense of the theory of Poisson-Nijenhuis structures. If
this happens, and if $AJ=\pm JA$, $M$ is foliated by ${\rm im}\,A$
into non degenerate K\"ahler-Nijenhuis submanifolds. If $A$ is a
non degenerate $(1,1)$-tensor field on $M$, $(M,g,J,A)$ is a
K\"ahler-Nijenhuis manifold iff one of the following two
properties holds: 1) $A$ is associated with a symplectic structure
of $M$ that defines a Poisson structure compatible with the
Poisson structure defined by $\Omega$; 2)  $A$ and $A^{-1}$ are
associated with closed $2$-forms. On a K\"ahler-Nijenhuis
manifold, if $A$ is non degenerate and $AJ=-JA$, $A$ must be a
parallel tensor field.}
\end{minipage}
\end{center} \vspace{5mm}
%\noindent
\section{Definition and basic formulas}
A K\"ahler manifold is a particular case of a symplectic
$2n$-dimensional manifold $(M,\Omega)$ with a symplectic form
defined as
\begin{equation} \label{Omega} \Omega(X,Y)=g(JX,Y)\hspace{5mm}(X,Y\in \Gamma TM),
\end{equation}
where $\Gamma$ denotes the space of global cross sections, $J$ is
a complex structure on $M$, and $g$ is a Hermitian metric on
$(M,J)$ \cite{KN}. Accordingly, on $M$ one has the Poisson
bivector field $\Pi$ defined by the Poisson brackets computed with
the symplectic form $\Omega$. The aim of this note is to discuss
Nijenhuis tensor fields $A$ that are compatible with $\Pi$ in the
sense of the theory of Poisson-Nijenhuis manifolds
\cite{{KSM},{MM},{V3},{V2}}. If this happens, the quadruple
$(M,g,J,A)$ will be called a {\it K\"ahler-Nijenhuis manifold},
and $A\in\Gamma\,{\rm End}(TM)$ will be a {\it K\"ahler-compatible
Nijenhuis} (K.c.N.) tensor field. The interest in
Poisson-Nijenhuis structures comes from their usefulness in the
search of first integrals of hamiltonian dynamical systems
\cite{MM}.

In what follows, we will use {\it musical morphisms} defined by
formulas of the type
\begin{equation} \label{musical} \beta(\sharp_\Pi\alpha)=
\Pi(\alpha,\beta),\;(\flat_\Omega X)(Y)=\Omega(X,Y),\end{equation}
where $\Pi$ may be any $2$-contravariant and $\Omega$ any
$2$-covariant tensor field. In particular, (\ref{Omega}) is
equivalent with $ \flat_g\circ J=\flat_\Omega$, and we have
$\sharp_\Pi\circ\flat_\Omega=-Id$ and
$\Pi=\sharp_\Pi\Omega=\sharp_g\Omega$ $(\sharp_g=\flat^{-1}_g)$
i.e., \begin{equation} \label{relPiOmega}
\Pi(\alpha,\beta)=\Omega(\sharp_\Pi\alpha,\sharp_\Pi\beta)=
\Omega(\sharp_g\alpha,\sharp_g\beta).
\end{equation}

Known results on Poisson-Nijenhuis structures \cite{{MM},{V2}}
tell that, on any symplectic manifold $(M,\Omega)$, the tensor
field $A\in\Gamma\,{\rm End}(TM)$ defines a Poisson-Nijenhuis
structure on $M$ iff $A=\sharp_\Pi\circ\flat_\Theta$, where
$\Theta$ is a closed, differential $2$-form such that one of the
following properties holds:\\ 1) $A$ is a Nijenhuis tensor field
i.e.,
\begin{equation} \label{NijA} {\rm Nij}_A(X,Y)=
[AX,AY]-A[X,AY]-A[AX,Y]+A^2[X,Y]=0;
\end{equation} 2) $\sharp_\Pi\Theta$ is a Poisson bivector field,
i.e., \begin{equation} \label{complementary}
[\sharp_\Pi\Theta,\sharp_\Pi\Theta]=0, \end{equation} where
$[\;,\;]$ is the Schouten-Nijenhuis bracket \cite{V1};\\ 3)
$\Theta$ is a complementary $2$-form of $\Pi$ i.e. \cite{V2},
\begin{equation} \label{comp2} \{\Theta,\Theta\}=0, \end{equation}
where $\{\;,\;\}$ is the Koszul bracket \cite{V1};\\ 4) the
$2$-form $\tilde\Theta$ defined by \begin{equation}
\label{tildeTheta}
\flat_{\tilde\Theta}=\flat_{\Theta}\circ\sharp_\Pi\circ\flat_\Theta
\end{equation} is closed.

Thus, if we add the request that $M$ is a K\"ahler manifold, the
above conditions characterize K.c.N. tensor fields. Furthermore, a
K\"ahler-Nijenhuis structure is also defined by the closed form
$\Theta$ with the properties 1)-4). We will say that $\Theta$ is
the {\it associated} K.c.N. form of $A$, which, in turn, is {\it
associated} with $\Theta$. Notice that
\begin{equation} \label{ThetaA} \Theta(X,Y)=-\Omega(AX,Y)=-\Omega(X,AY)
\end{equation} (use the skew symmetry of $\Theta$), and
\begin{equation} \label{tildeTheta2} \tilde\Theta(X,Y)=
-\Omega(AX,AY). \end{equation} In the rest of the paper, all the
encountered $(1,1)$-tensor fields $A$ are supposed to satisfy the
second equality (\ref{ThetaA}), called the $\Omega$-{\it
skew-symmetry} of $A$. $\Omega$-skew-symmetry ensures that
$A=\sharp_\Pi\circ\flat_\Theta$, where $\Theta$ is a $2$-form.

If $(x^i)$ $(i=1,...,2n)$ are local coordinates on $M$,
characteristic property 2) becomes (\cite{V1}, Proposition 1.5)
\begin{equation} \label{localexpr} \sum_{{\rm
Cycl}(i,j,k)}\Omega^{uv}\Theta_{ui}\nabla_v\Theta_{jk}
\stackrel{(\ref{ThetaA})}{=}\sum_{{\rm
Cycl}(i,j,k)}A_i^u\nabla_u\Theta_{jk}=0,
\end{equation} where $\nabla$ is the Levi-Civita connection of
$g$, and we use the Einstein summation convention. Thus, the
$\Omega$-skew-symmetric tensor field $A$ is K.c.N. iff
\begin{equation} \label{freeexpr} \sum_{{\rm Cycl}(X,Y,Z)}
(\nabla_{X}\Theta)(Y,Z)=0,\hspace{2mm}\sum_{{\rm Cycl}(X,Y,Z)}
(\nabla_{AX}\Theta)(Y,Z)=0, \end{equation} $\forall X,Y,Z\in
\Gamma TM$, where the first condition is equivalent to $d\Theta=0$
and the second condition is the coordinate-free equivalent of
(\ref{localexpr}). Notice also that, in view of (\ref{ThetaA}),
conditions (\ref{freeexpr}) are equivalent to
\begin{equation}\label{compatibilitate}
\sum_{Cycl(X,Y,Z)}\Omega[(\nabla_XA)(Y),Z]=0, \hspace{2mm}
\sum_{Cycl(X,Y,Z)}\Omega[(\nabla_{AX}A)(Y),Z]=0,\end{equation}
respectively, and the $\Omega$-skew-symmetric tensor field $A$ is
K.c.N. iff it satisfies (\ref{compatibilitate}).

On the other hand, characteristic property 3) has the interesting
equivalent form \cite{{V1},{V2}}
\begin{equation} \label{deltaC}
\delta^C(\Theta\wedge\Theta)=2(\delta^C\Theta)\wedge\Theta,
\end{equation} where $\delta$ is the Riemannian codifferential, and
$\delta^C=C\circ\delta\circ C$, with $C$ defined by the action of
$J$ on the arguments of a form (e.g., \cite{Gold}).

On a K\"ahler manifold, it is natural to consider the following
particular cases. We will say that a tensor field
$A\in\Gamma\,{\rm End}(TM)$ is {\it complex-compatible} (c.c.) if
$A\circ J=J\circ A$, and is {\it skew-complex-compatible} (s.c.c.)
if $A\circ J=-J\circ A$. Furthermore, if
$A=\sharp_\Pi\circ\flat_\Theta$ where $\Theta$ is a $2$-form, $A$
is c.c. iff $\Theta$ is of the complex type $(1,1)$ and $A$ is
s.c.c. iff $\Theta$ has components of the complex type $(2,0)$ and
$(0,2)$ only. This means that $\Theta(JX,JY)=\pm\Theta(X,Y)$,
respectively, and, if we denote by
\begin{equation} \label{projectors} \mathcal{P}=\frac{1}{2}(Id\otimes Id+J\otimes
J),\;\tilde{\mathcal{P}}=\frac{1}{2}(Id\otimes Id-J\otimes J)
\end{equation} the projectors of $2$-covariant tensors onto their
components of complex type $(1,1)$ and $[(2,0)+(0,2)]$ (each
factor of the tensor product acts on the corresponding argument),
such forms may be written as
\begin{equation} \label{proiectii} \Theta=\mathcal{P}\Xi,\,
\Theta=\tilde{\mathcal{P}}\Xi,\end{equation} respectively, where
$\Xi$ is an arbitrary $2$-form on $M$. In both cases, we are
speaking of a real form $\Theta$, and we will say that $\Theta$ is
c.c., in the first case, and s.c.c., in the second case.

In these two cases, the conditions that ensure the K.c.N. property
may be written under specific forms. Let us denote
\begin{equation} \label{defBC}\begin{array}{c}
E_A(X,Y)=(\nabla_XA)(Y),\,F_A(X,Y)=(\nabla_{AX}A)(Y),\vspace{2mm}\\
B_A=alt(E_A),\, C_A=alt(F_A), \end{array}\end{equation} where
$alt$ is the skew-symmetric part of a tensor. Then, we get
\begin{prop} \label{proposBC} 1. The $\Omega$-skew-symmetric,
c.c. tensor field $A$ is K.c.N. iff
\begin{equation} \label{cond1,1}\tilde{\mathcal{P}}B_A=0,
\,\tilde{\mathcal{P}}C_A=0.\end{equation} 2. The
$\Omega$-skew-symmetric, s.c.c. tensor field $A$ is K.c.N. iff
conditions {\rm(\ref{compatibilitate})} hold $\forall
X,Y,Z\in\Gamma T^cM$ $(T^cM=TM\otimes\mathbb{C})$ that are of the
complex type $(1,0)$, and
\begin{equation} \label{cond-1,1}
\mathcal{P}E_A=0,\,\mathcal{P}F_A=0.\end{equation}
\end{prop} \noindent{\bf Proof.} If $A$ is c.c. the extension of
$A$ to $T^cM$ preserves the complex type and, since $\nabla J=0$,
the same holds for the operators $\nabla_XA$, $\forall X\in\Gamma
T^cM$. Now, using the fact that $\Omega$ has the complex type
$(1,1)$, we see that conditions (\ref{compatibilitate}) are
identically satisfied if $X,Y,Z$ are of the same complex type.
Furthermore, from the $\Omega$-skew-symmetry of $A$ it follows
easily that, $\forall X\in \Gamma TM$, the tensor field
$\nabla_XA$ is also $\Omega$-skew-symmetric. This implies that,
for two arguments, say $X,Y$, of the same complex type (e.g.,
$(1,0)$) and the third, $Z$, of opposite type ($(0,1)$),
(\ref{compatibilitate}) is equivalent to $B_A(X,Y)=0, C_A(X,Y)=0$.
This happens iff (\ref{cond1,1}) holds.

Similarly, if $A$ is s.c.c., $A$ and $\nabla_XA$ change the
complex type of the vectors from $(1,0)$ to $(0,1)$ and
conversely. This implies that, for two arguments of the same type
and the third of the opposite type, (\ref{compatibilitate}) is
equivalent to $E_A(X,Y)=0,F_A(X,Y)=0$ whenever $X,Y$ have opposite
complex types. This is the same thing as conditions
(\ref{cond-1,1}). Of course, we must still ask
(\ref{compatibilitate}) to hold for three arguments of the same
complex type. Q.e.d.

Notice also that in the c.c. and s.c.c. cases (\ref{deltaC})
becomes
\begin{equation} \label{eqprop11}
\delta(\Theta\wedge\Theta)-2(\delta\Theta)\wedge\Theta=0.
\end{equation}

We end this section by a number of examples.
\begin{example} \label{example1} {\rm Any parallel $2$-form $\Theta$
of a K\"ahler manifold is a K.c.N. form. In particular, if a
K\"ahler manifold has a parallel Ricci tensor field, the Ricci
form is a c.c. form that defines a K\"ahler-Nijenhuis structure.}
\end{example}
\begin{example} \label{example2} {\rm Let $M$ be a
hyper-K\"ahler manifold with the metric $g$, the parallel complex
structures $(J_1,J_2,J_3)$ that satisfy the quaternionic
identities, and the respective K\"ahler forms
$\Omega_1,\Omega_2,\Omega_3$. Then, the tensors $J_2,J_3$ are
s.c.c., K.c.N. tensor fields on the K\"ahler manifold
$(M,g,J_1,\Omega_1)$. The corresponding K.c.N. forms are the
K\"ahler forms $-\Omega_3,\Omega_2$, which are parallel forms
\cite{Bg}.}
\end{example}
\begin{example} \label{example3}
{\rm On a compact Hermitian symmetric space, any real closed
$2$-form $\Theta$ which has no $(1,1)$-component is a K.c.N. form.
Indeed, $d\Theta=0$ implies that the $(2,0)$-component of $\Theta$
is holomorphic hence, harmonic (e.g., \cite{Gold}). Therefore,
$\Theta$ is harmonic, and, because the manifold is a compact
Hermitian symmetric space, $\Theta\wedge\Theta$ is harmonic too
(e.g., \cite{Helg}). Thus, $\Theta$ satisfies condition
(\ref{eqprop11}). Moreover, since $\Theta$ is s.c.c., by a result
that will be proven at the end of this paper, $\Theta$ is a
parallel form.}
\end{example}
\begin{example} \label{example4}
{\rm On a compact Hermitian symmetric space any real, harmonic
$(1,1)$-form $\Theta$  is a c.c., K.c.N. form. (Use again the
final argument of Example \ref{example3}).}
\end{example}
\begin{example}\label{example0} {\rm On $M=\mathbb{C}^n$, with the
flat K\"ahler metric and the natural complex coordinates
$(z^\alpha)$, the $(1,1)$-form $\Theta=z^1dz^1\wedge d\bar z^2$ is
closed and satisfies condition (\ref{eqprop11}). Hence, $\Theta$
is a c.c., K.c.N. form. It is easy to check that $\Theta$ is not a
parallel form.}
\end{example}
\section{Geometric Properties}
Let $(M,g,J,A)$ be a K\"ahler-Nijenhuis manifold. The basic
geometric object that we detect beyond the usual K\"ahlerian
objects is the differentiable, generalized distribution $
\mathcal{A}={\rm im}\,A$. It is well known that this distribution
is completely integrable. For the record, we write down a
straightforward proof below.
\begin{prop} \label{prop21} If $A$ is a Nijenhuis tensor field
(i.e., {\rm(\ref{NijA})} holds), the generalized distribution $
\mathcal{A}={\rm im}\,A$ is completely integrable. \end{prop}
\noindent{\bf Proof.} Condition (\ref{NijA}) shows that $
\mathcal{A}$ is an involutive distribution. Hence, integrability
will follow from the Sussmann-Stefan-Frobenius theorem (e.g.,
\cite{V1}) if we prove that $ \mathcal{A}$ is invariant i.e.,
$\forall X,Y\in\Gamma TM$ one has
$[\exp(tAX)]_*(AY)\in\mathcal{A}$, $\forall t\in{\bf R}$ such that
$\exp(tAX)$ exists.

Denote $A_x(t)=[\exp(tAX)]_*(A_{\exp(-tAX)(x)})$ $(x\in M)$. Then,
\begin{equation} \label{difeq1} \frac{dA_x(t)}{dt}=
\lim_{s\rightarrow0}\frac{1}{s}\{[\exp((t+s)AX)]_*(A_{\exp[-(t+s)AX](x)})
\end{equation}
$$-[\exp(tAX)]_*(A_{\exp(-tAX)(x)})\}=[L_{AX}A(t)]_x,$$
where $L$ denotes the Lie derivative.

The required invariance of $ \mathcal{A}$ will be a consequence of
the local existence of a $(1,1)$-tensor field $C(t)$ such that
$A(t)=A\circ C(t)$. If $C(t)$ exists, (\ref{difeq1}) implies $$
A\circ\frac{\partial C}{\partial t}=L_{AX}(A\circ
C)=(L_{AX}A)\circ C+A\circ L_{AX}C $$
$$\stackrel{(\ref{NijA})}{=}A\circ L_XA\circ C+A\circ L_{AX}C.$$
Therefore, if $C(t)$ satisfies
\begin{equation} \label{eqdif}
\frac{\partial C}{\partial t}-(L_XA)\circ
C-L_{AX}C=0,\;C(0)=Id,\end{equation} $\forall x\in M$, $A(t)$ and
$A\circ C(t)$ satisfy the same differential equation
(\ref{difeq1}) and the same initial condition, and must be equal.
Since (\ref{eqdif}) has a local solution, we are done. Q.e.d.

Thus, through every point $x\in M$ one has a {\it characteristic
leaf}, the maximal integral submanifold of the generalized
distribution $ \mathcal{A}$, which we denote by $
\mathcal{L}=\mathcal{L}_x$, immersed in $M$ by
$\iota=\iota_\mathcal{L}:\mathcal{L}\hookrightarrow M$.
\begin{prop} \label{thfoilor} Let $(M,g,J,A)$ be a c.c. or
s.c.c. K\"ahler-Nijenhuis manifold. Then, each characteristic leaf
$ \mathcal{L}$ inherits an induced structure of a non degenerate
K\"ahler-Nijenhuis manifold with the normal bundle $
\mathcal{K}|_\mathcal{L}$, where $ \mathcal{K}= {\rm ker}\,A={\rm
ker}\,\Theta$. Furthermore, if the structure $A$ is regular, the
decomposition $TM=\mathcal{A}\oplus \mathcal{K}$ is a complex
analytic, orthogonal, locally product structure on $M$.
\end{prop} \noindent{\bf Proof.}
In the case of a c.c. or s.c.c. tensor field $A$, the distribution
$ \mathcal{A}$ is $J$-invariant, and the characteristic leaves $
\mathcal{L}$ are K\"ahler submanifolds of $M$. Then, $\forall
X,Y\in\Gamma TM$, we have
\begin{equation}\label{relgTheta}
g(AX,Y)\stackrel{(\ref{Omega})}{=}-\Omega(JAX,Y) =
\mp\Omega(AJX,Y) \stackrel{(\ref{ThetaA})}{=}\mp\Theta(Y,JX),
\end{equation} and we see that $Y\in\mathcal{K}={\rm ker}\,\Theta$ iff
$Y\perp\mathcal{A}$. Therefore, the normal bundle of $
\mathcal{L}$ is $ \mathcal{K}|_\mathcal{L}$. Since
$A=\sharp_\Pi\circ\flat_\Theta$ and $\sharp_\Pi$ is an
isomorphism, we also have ${\rm ker}\,A={\rm ker}\,\Theta$.

The field of planes $ \mathcal{K}$, which, by the above result, is
$J$-invariant, is not a differentiable distribution since its
dimension is not lower semi-continuous. Differentiability occurs
iff the c.c. or s.c.c. Nijenhuis tensor $A$ (and the corresponding
form $\Theta$) is {\it regular} i.e., of a constant rank. Then, as
it is well known, $d\Theta=0$ implies that $ \mathcal{K}$ is
involutive, and the decomposition
$TM=\mathcal{A}\oplus\mathcal{K}$ defines a complex analytic,
orthogonal, locally product structure on $M$.

Of course, the distribution $ \mathcal{A}$ also is invariant by
$A$ hence, $A|_\mathcal{L}$ is a $(1,1)$-tensor field on $
\mathcal{L}$. Moreover $A|_\mathcal{L}$ is a Nijenhuis tensor,
since the Lie brackets of condition (\ref{NijA}) are
$\iota$-compatible. In the c.c. and s.c.c. cases $A|_\mathcal{L}$
has zero kernel because ${\rm ker}\,A\perp T\mathcal{L}$. Hence,
$A|_\mathcal{L}$ is non degenerate, and so is the associated form
$\Theta_\mathcal{L}$. Furthermore, formula (\ref{ThetaA}) shows
that $\Theta_\mathcal{L}=\iota^*\Theta$, and property 1) of the
K.c.N. structures shows that $A|_\mathcal{L}$ is a K.c.N. tensor
field. Q.e.d.

If the Nijenhuis tensor $A$ is non degenerate, the manifold $M$
itself is the only characteristic leaf. Furthermore, if
$(M,g,J,A)$ is a non degenerate K\"ahler-Nijenhuis manifold the
corresponding K.c.N. form $\Theta$ is a symplectic form and the
Poisson brackets of the latter define a Poisson bivector field
$\Psi$.
\begin{prop} \label{proposcompat} Let $A$ be a $\Omega$-skew-symmetric,
non degenerate $(1,1)$-tensor field on $M$. Then: 1. $A$ is K.c.N.
iff it is associated with a closed $2$-form $\Theta$ and the
Poisson structures defined by $\Omega,\Theta$ are compatible i.e.,
$[\Pi,\Psi]=0$.\\ 2. $A$ is K.c.N. iff both $A$ and $A^{-1}$ are
associated with closed $2$-forms.
\end{prop} \noindent{\bf Proof.} For {\it 1}, by a result proven in \cite{MM}
(see also \cite{V3}), the compatibility condition $[\Pi,\Psi]=0$
implies the fact that $A$ is K.c.N. Conversely, from
$A=\sharp_\Pi\circ\flat_\Theta$ we get
\begin{equation}\label{eqA-1}\sharp_\Psi=A^{-1}\circ\sharp_\Pi.
\end{equation} Hence
the Poisson structure $\Psi$ belongs to the enlarged Poisson
hierarchy of the Poisson-Nijenhuis structure $(\Pi,A)$, and the
required compatibility follows from the properties of the Poisson
hierarchy (e.g., \cite{V3}).

For {\it 2}, again, the theorem of the Poisson hierarchy tells us
that if $A$ is K.c.N. then $A^{-1}$ is K.c.N. too. Hence, if we
write $A^{-1}=\sharp_\Pi\circ\flat_{\Theta'}$, $\Theta'$ must be
closed. Conversely, assume that $\Theta$ and $\Theta'$ are closed.
From (\ref{eqA-1}), we get
\begin{equation} \label{eqauxil}
\flat_{\Theta'}=\flat_\Omega\circ\sharp_\Psi\circ\flat_\Omega,
\end{equation} therefore, by property 4) of the K.c.N. structures
(see Section 1), $(\Psi,A^{-1})$ is a Poisson-Nijenhuis structure,
and $\Pi$ belongs to the Poisson hierarchy of the former.
Therefore, $[\Psi,\Pi]=0$ and, by part {\it 1} of the proposition,
we are done. (In fact, since $\Theta$ is a symplectic form, and in
view of (\ref{eqA-1}), $(\Psi,A^{-1})$ is a Poisson-Nijenhuis
structure iff the Poisson bivector fields $\Psi,\Pi$ are
compatible.) Q.e.d.
\begin{rem} \label{obspt1,1} {\rm In the c.c. case, conclusion {\it
2} of Proposition \ref{proposcompat} follows immediately from the
first part of Proposition \ref{proposBC}. Indeed, we can use
$\nabla(A\circ A^{-1})=0$ to derive
\begin{equation}\label{1,1BC} C_A(X,Y)=
A(B_{A^{-1}}(AX,AY)),\end{equation} and conclude as required from
(\ref{cond1,1}). We also notice the formula
\begin{equation} \label{NijcuB}
-A^{-1}({\rm Nij}_A(X,Y))=2
[B_{A^{-1}}(AX,AY)+B_A(X,Y)].\end{equation}} \end{rem}
\begin{corollary} \label{corolfinal} Let $A\in {\rm End}(TM)$
define a c.c., orthogonal, almost product structure on $M$. Then,
the tensor field $A$ is associated with a $(1,1)$-form $\Theta$,
and $A$ is K.c.N. iff $\Theta$ is closed
\end{corollary} \noindent{\bf Proof.} The orthogonality of the
structure means $g(AX,AY)=g(X,Y)$, and it implies that
$\Theta(X,Y)=-\Omega(AX,Y)$ is skew symmetric. Thus, $\Theta$ is
the required $2$-form. Furthermore, $A^{-1}=A$, and the result
follows from part {\it 2} of Proposition \ref{proposcompat}.
Q.e.d.

In particular, if $\Theta$ of Corollary \ref{corolfinal} is closed
$Nij_A=0$ and the almost product structure $A$ is integrable.

We finish by showing that for the non degenerate, s.c.c. tensor
fields the K.c.N. condition is very restrictive.
\begin{prop} \label{theorem21}
A non degenerate, $\Omega$-skew-symmetric, s.c.c. tensor field
$A\in\Gamma\,{\rm End}(TM)$ is K.c.N. iff $A$ is parallel.
\end{prop} \noindent{\bf Proof.} The quickest way to conclude is
by a local computation. Consider local, complex analytic
coordinates $(z^\alpha)$ $(\alpha=1,...,n)$ on $M$. The s.c.c.
property of $A$ means that the only possibly non-zero components
of $A$ are $(A^{\bar\beta}_\alpha,\,A^\beta_{\bar\alpha})$, and
$\Theta$ has no component of the complex type $(1,1)$. Since
$d\Theta=0$, the complex $(2,0)$-component of $\Theta$ is
holomorphic. Accordingly, condition (\ref{localexpr}) becomes
\begin{equation} \label{analy-anti} A^\alpha_{\bar\lambda}
\nabla_\alpha\Theta_{\mu\nu}=0,
\end{equation} and if $A$ is non degenerate we get
$\nabla_\alpha\Theta_{\mu\nu}=0$. Q.e.d.
\begin{rem}\label{ultima} {\rm Except for Proposition
\ref{thfoilor}, the results of this note also hold for {\it
pseudo-K\"ahler manifolds} i.e., where the metric $g$ is non
degenerate but it may not be positive definite.} \end{rem} \hspace*{7.5cm}{\footnotesize \begin{tabular}{l} Department of Mathematics\\
University of Haifa, Israel\\ E-mail: vaisman@math.haifa.ac.il
\end{tabular}}
\end{document}